\documentclass[journal]{IEEEtran}
\pdfoutput=1
\usepackage{framed} 
\usepackage{color}
\usepackage{graphicx}
\usepackage{subfigure}
\usepackage{mathtools}
\usepackage{paralist}
\usepackage{cleveref}
\usepackage{amsmath}
\newcommand {\eqrefn}{Eq.~\eqref}
\begin{document}
\title{Integrated bio-electrochemical model for a micro photosynthetic power cell}
\author{Tanneru~Hemanth~Kumar,~Resmi~Suresh~M.P,~Aravind~Vyas~Ramanan,~\IEEEmembership{Student~Member,~IEEE},~Shahparnia.M,~\IEEEmembership{Student~Member,~IEEE},~Muthukumaran~Packirisamy,~\IEEEmembership{Member,~IEEE},~Pragasen~Pillay~\IEEEmembership{Fellow,~IEEE},~Sheldon~Williamson,~\IEEEmembership{Senior~Member,~IEEE},~Philippe~Juneau and~Raghunathan~Rengaswamy~\thanks{Tanneru Hemanth Kumar, Resmi Suresh ~M.P and Raghunathan Rengaswamy is with SENAI, Department of Chemical Engineering, Indian Institute of Technology Madras, Chennai, Tamil Nadu, India 600036.}
\thanks{Aravind Vyas Ramanan,Shahparnia.M, and MuthukumaranPackirisamy is with Optical Bio-MEMS laboratory, Department of Mechanical Engineering, Concordia University, Montr$\acute{e}$al, QC, H3G1M8, Canada.} 
\thanks{Pragasen~Pillay and Sheldon~Williamson, is with P.D Ziogas Power Electronics laboratory, Department of Electrical and Computer Engineering, Concordia University, Montr$\acute{e}$al, QC, H3G1M8, Canada} 
\thanks{Philippe~Juneau is with Department of Biological Sciences, Universit$\acute{e}$ du Qu$\acute{e}$bec $\grave{a}$ Montr$\acute{e}$al (UQAM), Montr$\acute{e}$al, QC, H3C 3P8, Canada.}}

\maketitle

\begin{abstract}
A simple first-principles mathematical model is developed to predict the performance of a micro photosynthetic  power cell ($\mu$PSC), an electrochemical device which generates electricity by harnessing electrons from photosynthesis in the presence of light. A lumped parameter approach is used to develop a model in which the electrochemical kinetic  rate constants and diffusion effects are lumped into a single characteristic rate constant $K$. A non-parametric estimation of $K$ for the $\mu$PSC is performed by minimizing the sum square errors (SSE) between the experimental and model predicted current and voltages. The developed model is validated by  comparing the model predicted  $v-i$ characteristics with experimental data not used in the parameter estimation. Sensitivity analysis of the design parameters and the operational parameters reveal interesting insights for performance enhancement. Analysis of the model also suggests that there are two different operating regimes that are observed in this $\mu$PSC. This modeling approach can be used in other designs of $\mu$PSCs for performance enhancement studies.
\end{abstract}


\section{Introduction}

\IEEEPARstart{E}{nergy} consumption is increasing all over the world. Fossil fuels being being a major source of energy are getting depleted much faster than they can be replenished. As a result, renewable energy sources are being researched for applications that require different power ranges. Low power range application devices such as remote location sensors, bio-sensors are of much interest in recent years. Scaling down high power range devices for such applications is difficult due to various design issues. Microbial fuel cells, which are electrochemical devices that use electrons produced during respiration of microbes to generate current \cite{siu2008microfabricated, cheng2011increasing} has gained much attention for such applications. Mathematical modeling  of microbial fuel cells have also been attempted\cite{zhang1995modelling,picioreanu2010modelling,pinto2011unified}.

Micro photosynthetic cells ($\mu  PSC$) are a sustainable option for low power applications. $\mu  PSC$  uses  oxygenic photosynthetic organisms such as algae to generate current in the presence of light and function as a microbial fuel cell in the absence of light. A major advantage of  $\mu  PSC$ is its ability to generate current by harnessing the electrons from the electron transport chains in photosynthetic organelle of photoautotrophs using sunlight. In the absence of light, the cell generates current from the electrons that are harnessed from the metabolic pathways of the respiration process in photosynthetic organisms. Prototypes for $\mu  PSC$ \cite{yagishita1996photosynthetic,lam2004bio,rosenbaum2005utilizing,lam2006mems,chiao2006micromachined, shahparnia2011polymer, arvind2014advancedfabric} are available in the literature.

Till date, the focus has been mainly on experimental aspects of $\mu PSC$, recent works reported on $\mu PSC$ \cite{ramanan2015advanced},\cite{shahparnia2015micro}, and there has been very little attempt at developing mathematical models for $\mu PSC$. In general, the aim of any modeling exercise is to understand the underlying physical phenomena of a device and explore methods for improving device performance. Modeling of  $\mu PSC$ can help in understanding the performance limiting step(s) in a series of processes that occur during the operation of the device. Performance enhancement of the device can be achieved by focusing on the rate limiting steps. Modeling also helps in determining the optimal design and operational parameters, that can maximize the device performance. 

Modeling a system like $\mu PSC$ is complex since the device performance depends on the interactions of microorganisms with the operational parameters such as light intensity, quantum yield and so on. Further,  design parameters such as electrode structure and  the electrochemical interactions at the surface of the electrodes have an effect on the device performance. A mathematical model which incorporate all the phenomena that occur during the operation of $\mu PSC$ will be complicated. 

In this work, a simple model based on first principles is developed. The aim of this modeling exercise is to predict the  performance of $\mu PSC$, given device specifications and a set of operational parameters. 
In the present work, a lumped parameter model approach is used in the model development. The number of parameters chosen to describe the various processes will be largely determined by the richness of data in terms of the variables that are measured. The lumped parameter used in the current study is the characteristic rate constant $K$. 

This paper is organized as follows. The section on operation principles describes the details of $\mu PSC$ working. The model equations are described in the model formulation section. This is followed by a description of the methodology adopted for solving the model equations and parameter estimation from the published \cite{shahparnia2011polymer} $v-i$ data. Subsequently, model validation and sensitivity analysis studies that are performed are described. Finally, the utility of the model and interesting insights that can be derived from such a model analysis are outlined.

\section{Operation principle}
\label{experiments}
\subsection{Operation of {$\mu PSC$}}
\begin{figure}
\centering
\includegraphics[scale=0.5]{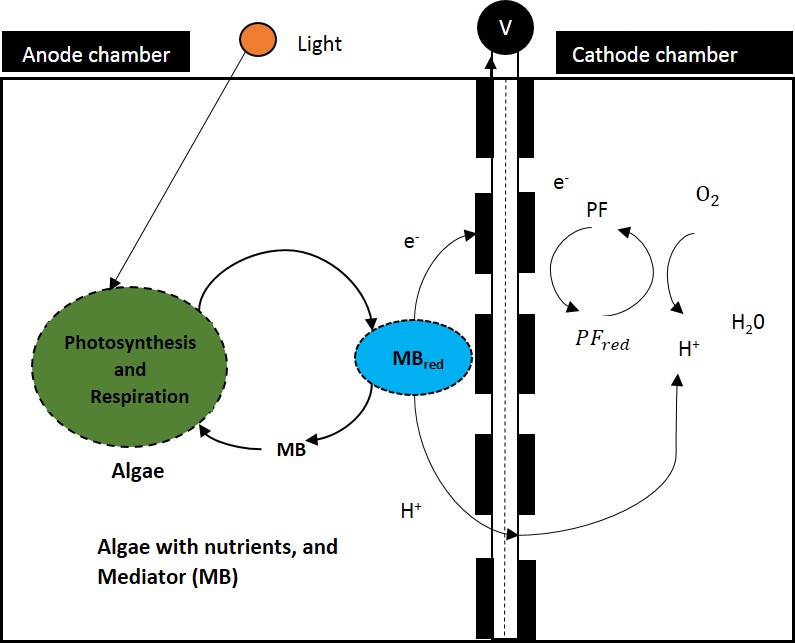}
\caption{Working of Photosynthetic cell: Electrons and protons released during photosynthesis and respiration are siphoned by the mediator to the electrode surface. Electrons flow to the cathode chamber through the external circuit to produce electricity. Protons diffuse through the membrane and reacts with the electrons and oxygen to produce water in the cathode chamber}
\label{Schematic}
\end{figure}

A schematic of the photosynthetic  cell is shown in Figure \ref{Schematic}. A $\mu PSC$ device consists of two chambers (anode and cathode) separated by a proton exchange membrane (Nafion). Anode and cathode chambers have a capacity to hold $2 ml$ of anolyte and catholyte respectively. Porous gold electrode patterns of $100 \mu\!m$ thick developed on both the sides of the Nafion membrane using lithography techniques act as both electrodes and current collectors.

Green algae (\textit{Chlamydomonas reinhardtii}) suspended in a growth medium (Sueoka's high salt medium, HSM) and a mediator (methylene blue) are the major components of the anolyte. Potassium ferricyanide ($PF$) solution is used as catholyte. Cell growth (using both nutrients and glucose as substrate) and cell decay occur inside the anode chamber.

Photosynthesis takes place in the presence of light, producing glucose from carbon dioxide and water. Respiration occurs in both dark and light conditions. The reactions are as follows.

\begin{equation} \label{eq.photo}
Photosynthesis: 6CO_2 + 6H_2O \rightarrow C_6H_{12}O_6 + 6O_2
\end{equation}

\begin{equation} \label{eq.respi}
Respiration: C_6H_{12}O_6 + 6O_2 \rightarrow 6CO_2 + 6H_2O 
\end{equation}

Electrons and protons are released during both respiration and photosynthesis. The mediator methylene blue (MB) diffuses into the microorganism and siphons these electrons from the electron carriers NADPH (during photosynthesis) / NADH (during respiration) during which it gets reduced (see \eqrefn{eq.re}) The reduced methylene blue ($MB_{red}$) then diffuses out of the microorganism and releases the electrons at the anode surface along with the protons thereby converting back to its original oxidized form (see \eqrefn{eq.res}). The electrons from the anode travel through the external circuit to the cathode chamber  producing electricity. Protons diffuse through the Nafion membrane to the cathode side.
\begin{multline}\label{eq.re}
MB + \hspace{0.01in}NADPH(photosynthesis)/NADH(respiration)\\
+ \hspace{0.01in}H^+   \rightarrow  MB_{red}+ NADP^+ / NAD^+
\end{multline}

\begin{equation}\label{eq.res}
 MB_{red} \rightarrow MB + 2e^- + 2H^+
\end{equation}

At the cathode surface, potassium ferricyanide($PF$) gets reduced to potassium ferrocyanide ($PF_{red}$) using the electrons from the external circuit, (see \eqrefn{fe3}). $PF_{red}$ thus formed supplies electron to $O_2$ and protons which combine to form water and $PF$ in the cathode chamber (see \eqrefn{fe2}).

\begin{equation} \label{fe3}
2PF + 2e^- \rightarrow 2PF_{red}
\end{equation} 

\begin{equation} \label{fe2}
2PF_{red} + 2H^{+} + 0.5O_2 \rightarrow 2PF + H_2O 
\end{equation}
In this paper, data from the experimental results reported in the prior work of the some of the authors of this paper \cite{shahparnia2011polymer,arvind2014advancedfabric,ramanan2015advanced,shahparnia2015micro} are used for parameter estimation and model validation.  

\section {Model formulation} \label{model formulation}

\subsection{Model equations}
\label{model}
The mathematical model developed in this work is intended to predict the current-voltage behavior of a $ \mu PSC $ in the presence of light. The inputs to the model are the initial concentrations of the species in the anode, light intensity, the external loads and the design parameters of the device.\\
\textit{Assumptions:}
\begin{enumerate}[(1)]
\item $\mu PSC$ is operated under isothermal conditions at $ 25^{0}C $ and atmospheric pressure.
\item Both anode and cathode chambers are assumed to be well mixed batch reactors.
\item Cell growth in anode chamber is governed by Monod kinetics. (Monod kinetics is most generally used mathematical model to  describe the growth of the suspended microorganisms in the aqueous medium as a function of concentration of nutrient medium)
\item Electrode kinetics follows Butler-Volmer equation.
\item Diffusion of  carbon dioxide, oxygen and sugar through Nafion membrane are assumed to be negligible.
\item Photosynthesis is considered to be the dominant process in the presence of light.
\item Concentration of species on the electrode surface is assumed to be equal to their corresponding bulk concentrations. No diffusion effects are considered.
\item Water is assumed to be in excess in the anode chamber.
\item Oxygen is assumed to be available in excess in the cathode chamber.
\item Activation losses near the cathode are assumed to be negligible.
\end{enumerate}
\subsubsection{Anode chamber}\label{subsec.anode}
The following phenomena occur in the bulk of the anode chamber. Cells grow by consuming the nutrient medium and they decay at a specified rate or when the nutrient medium is exhausted. With the assumptions stated, the temporal variation of concentrations of the species $x$ (cells), $N$ (nutrients), in the anode chamber can be described by the following differential equations.

The change in cell concentration $x$  by growth and decay of cells can be described by ~\eqrefn{cell}
\begin{equation} \label{cell}
\frac{d x}{dt} = k_{1}x - k_{2}x
\end{equation}
with $k_{2}$, the death rate of cells, $k_{1}$, the growth rate of cells, characterized by Monod kinetics.
\begin{equation*}
k_{1} = \mu_{max}\left(\frac{N}{K_{N}+N}\right)
\end{equation*}
where {$\mu_{max}$} is maximum specific growth rate, $N$ the nutrient medium used for growth, and $K_{N}$ is the half saturation constant w.r.t $N$.

The nutrient concentration change can be represented by ~\eqrefn{nutrient}
\begin{equation} \label{nutrient}
\frac{d N}{dt} = -\frac{k_{1}x}{Y_{x,N}}
\end{equation} 
where $Y_{x,N}$ is yield coefficient of cells w.r.t nutrients.

The next step is to understand the source of electrons and the mechanism by which they reach the anode surface. There are several reactions that occur during photosynthesis and a complete description of the process with detailed mechanisms can be very complicated. Therefore, we propose a simple one-step mechanism that results in a tractable and useful model. A good description of processes that happen during photosynthesis can be found in \cite{taiz2010plant}.
 
The first step in photosynthesis is the water splitting reaction. This reaction happens inside the thylakoid membrane of chloroplast. 
\begin{equation} \label{water split}
H_2O  \rightarrow  2H^{+} + 2e^- + 0.5 O_2  
\end{equation} 
The electrons are received by the electron acceptor $NADP^+ $ and forms $ NADPH $ at the end of the electron transport chain. Each $ NADP^+ $ can take $ 2e^- $ and one $H^+$.
\begin{equation} \label{NADPH formation}
NADP^+ + 2H^{+} + 2e^- \rightarrow  NADPH+ H^+    
\end{equation} 
MB added in the anolyte diffuses in to the cell and siphons the electrons and the protons from  $NADPH + H^+$ to form $MB_{red}$.
\begin{equation} \label{MB_red formation}
MB + \hspace{0.01in} NADPH+ H^+\rightarrow MB_{red}+ NADP^+
\end{equation}
At the anode surface $ MB_{red} $ is oxidized back to $MB$. 
\begin{equation} \label{MB_red to MB}
 MB_{red} \rightarrow MB+2H^++2e^-
\end{equation} 
The ideal way to model this phenomena is to consider the concentration variations of all the species in anode chamber, the rate of reactions occurring in bulk and the effect of diffusion on the concentration of species at the electrode surface for use in Bulter-Volmer equation.

To simplify the model, we conceptualize the electron carrier, $ MB $ and $ MB_{red}$ as intermediates and adding \eqrefn{water split}, \eqrefn{NADPH formation}, \eqrefn{MB_red formation} and \eqrefn{MB_red to MB}, we obtain 
\begin{equation}\label{Anode final}
H_2O  \rightarrow  2H^{+} + 2e^- + 0.5 O_2  
\end{equation}
Assuming \eqrefn{Anode final} occurs at the anode surface, the complexity will be reduced to a great extent by considering the cell concentration and light intensity as reactants. This is because, at a macro level, the rate at which the water splitting occurs in a microorganism is a function of both cell concentration and light intensity. An important point to note at this juncture is that the water splitting reaction is assumed as a representative of all the phenomena that occur in the anode chamber. This can be summarized in \eqrefn{Anode complete reaction}
\begin{equation} \label{Anode complete reaction}
x(cells) + H_2O\xrightarrow{Light,k_a}{}x(cells)+2H^++2e^-+0.5O_2 
\end{equation} 
Here, $ k_a$ has much more significance than a mere rate constant of the water splitting reaction since this rate constant represents the reaction rate and also the many transport phenomena that are involved in the movement of all the active species.
\subsubsection{Cathode chamber}\label{subsec.cathode}
The following processes occur in the cathode chamber. The electrons received at the cathode surface are used to reduce potassium ferricyanide($PF$) to potaassium ferrocyanode($ PF_{red} $). 
\begin{equation} \label{PFred}
2PF + 2e^- \rightarrow 2PF_{red}
\end{equation} 
Protons diffuse from the anode side to cathode chamber through Nafion, and take part in the reaction, where $PF_{red}$ is oxidized to $PF$ by donating electrons to oxygen to form water.  

\begin{equation} \label{PF}
2PF_{red} + 2H^{+} + 0.5O_2 \rightarrow 2PF + H_2O 
\end{equation}

Detailed modeling of the cathode chamber should ideally track the concentrations of all the species in cathode chamber and the influence of diffusion effects on the species concentrations at the electrode surface where the electrochemical reactions occur. This complexity can be handled if we assume that the overall reaction that occurs on the cathode surface as the oxygen reduction reaction (ORR). The rationale for this assumption is the same as the one used in the modeling of the anode chamber. 

Adding \eqrefn{PFred}  and \eqrefn{PF} we obtain 
\begin{equation}  \label{cathodefinal}
0.5O_2 + 2H^{+}+2e^- \xrightarrow{k_c} H_2O
\end{equation}
Similar to the anode chamber, the rate constant $ k_c $ has to be interpreted as not just being the rate constant for the oxygen reduction reaction.

\subsubsection{Electrochemical equations}
The voltage, $v$ and current, $i_{\mu PSC}$ produced when an external resistance, $R_{ext}$ is connected to a $\mu PSC$ is given by 
\begin{equation}\label{voltage}
v = i_{\mu PSC}R_{ext} = E_0 - \eta_{ohm}-\eta_{conc}-\eta_{act}
\end{equation}
 where $E_{0}$ is non-standard thermodynamic voltage;  $\eta_{ohm},\  \eta_{conc}$ and $\eta_{act}$  are ohmic, concentration and activation losses across the cell respectively.
Ohmic losses occur due to the transfer of current through the internal resistance of the $\mu PSC$ and can be represented by Ohm's law.
\begin{equation}\label{ohmic}
\eta_{ohm} =  i_{\mu PSC}R_{int}
\end{equation}

Using \eqrefn{ohmic}, and rewriting the activation losses at two electrodes individually, \eqrefn{voltage} takes the form  
\begin{equation}\label{v1}
i_{\mu PSC}R_{ext} = (E_C -\eta_{C,act})-(E_A +\eta_{A,act})- i_{\mu PSC}R_{int}-\eta_{conc}
\end{equation}
Concentration over-potentials are due to mass transport losses at higher current densities. Since the diffusion effects are not explicitly modeled, the concentration losses are incorporated into the characteristic rate constant and are not represented in terms of voltage.
 
Activation losses near cathode are neglected following \cite{noren2005clarifying} and  rewriting \eqrefn{v1} for  $\eta_{A,act}$, the anodic activation loss is related to current density of the anodic reaction.
 
\begin{equation}\label{volt2}
\eta_{A,act} = (E_C -E_A) -i_{\mu PSC}(R_{ext}+R_{int})
\end{equation}
with $E_C -E_A = E_0 $, \eqrefn{volt2} can be written as

\begin{equation}\label{volt3}
\eta_{A,act} = E_0 -i_{\mu PSC}(R_{ext}+R_{int})
\end{equation}

The current density  at the anode surface is given by the Butler-Volmer equation \cite{sunden2005transport}. 

\begin{multline} \label{BV}
j = \frac{i_{\mu PSC}}{A_{E}}\\ = j_{0}\left[exp\left(\frac{\alpha nF}{RT}\eta_{act,a}\right)-exp\left(\frac{-(1-\alpha)nF}{RT}\eta_{act,a}\right)\right] 
\end{multline} 
where $j$ is current density; $j_{0}$ is exchange current density and is given by \eqrefn{currdens} and $A_E$ is the electrode surface area.
 
\begin{equation} \label{currdens}
j_0 =k_a^\alpha~k_c^{(1-\alpha)}C_a^{\alpha}C_c^{(1-\alpha)}~= KnF\left(\frac{L_0~C_f~A_s~Q~\eta_{eff}~x}{N_{Av}~x_{max}}\right)^\alpha 
\end{equation}
with $K (\frac{1}{m^2}) $,~the characteristic rate constant,~$ L_0 $,~Intensity of incident light ($lux$), ~$C_f $, conversion factor (lumen to photons per sec),~$Q$,~quantum yield,($\frac{\text{number of electrons~ released}}{\text{number of photons absorbed}}$),~$k_a$, rate constant of anode reaction ,~$k_c$, rate constant of cathode reaction,~$C_a$, product of reactant species concentration in  anode reaction,~$C_c$, product of reactant species concentration in  cathode reaction and 
$A_s$, exposure surface area.~$C_c$ is taken as unity based on assumption that oxygen is in excess in cathode chamber.

Solving \eqrefn{volt3} and \eqrefn{BV} simultaneously, the current density of  $\mu PSC$ can be obtained. Voltage from the $\mu PSC$ can be calculated from \eqrefn{voltage}.
The characteristic rate constant, $K$, has information about both the rate constants $k_a$ and ${k_c}$ of the proposed consolidated equations for anode and cathode.  $K$, $k_{a}$ and $k_b$ are related by the equation $K=k_a^\alpha~k_c^{(1-\alpha)}$.

\subsection*{$E_{0} \hspace*{.2cm}calculation$:}The following reactions occur on the electrode surfaces. The standard reduction potentials of the bio-reactions at the anode are adapted from \cite{nelson2008lehninger,arvia2011electrochemical}.

Anode surface: Oxidation of reduced methylene blue to methylene blue. 
\begin{equation*}
MB_{red} \rightarrow MB + 2e^{-}+2H^{+}  \hspace*{1cm} E_{ox}^{0}  = -0.011 V
\end{equation*}

Cathode surface: Reduction of potassium ferricyanide to potassium ferrocyanide.
\begin{align*}
2PF +2e^{-} \rightarrow   2{PF_{red}} \hspace{0.5cm}  E_{red}^{0} &= 0.361 V \\
2{PF_{red}}+0.5O_{2}+2H^{+} \rightarrow 2PF + H_{2}O\\
2{PF_{red}} \rightarrow   2PF +2e^{-}  \hspace{0.5cm} E_{red}^{0} &= -0.361 V\\
0.5O_{2}+2H^{+} +2e^{-} \rightarrow  {H_{2}O}  \hspace{0.5cm}E_{red}^{0}&= 1.23 V\\
\end{align*}
\begin{equation*}
E_{cell}^{0} = E_{cathode}^{0}-E_{anode}^{0}= 1.23-(-0.011)= 1.241 V
\end{equation*}
Following \cite{gunawardena2008performance}, $E_{cell}^{0}$, the standard cell potential can be obtained at $STP$ with species concentrations at $1$M.  Generally, Nernst equation is used to relate the standard cell potential and the potential that can be obtained with the cell operating conditions. In the present study, the model developed being lumped, the effect of species concentrations at the operating conditions cannot be incorporated through the Nernst equation. Hence, the standard cell potential is used in the simulations assuming that the other terms are absorbed in the characteristic rate constant $K$.
\section{Solution to model equations}\label{solution}

The model equations presented contains $2$ ODEs (\eqrefn{cell} and \eqrefn{nutrient}) and $2$ algebraic equations (\eqrefn{volt3} and \eqrefn{BV}). The four unknown variables are: the cell concentration $x$, nutrient concentration $N$, the current density $j$, and the activation over potential $\eta_{a}$.

First, the two ODEs are integrated using the MATLAB inbuilt integrator (ODE15S) and the final concentration  $x$ at each iteration of $R_{ext}$ is used as an initial guess for the next  $R_{ext}$ and also to calculate $j_0$. The current density of $\mu PSC$,~$j$, is obtained by solving the two algebraic equations using the non-linear equation solver of MATLAB. The current and the voltage from the model are calculated by using \eqrefn{BV} and \eqrefn{voltage} respectively.
 
\section{Parameter estimation} \label{parameterest}
Various parameters present in the model and their approximate values taken from the literature are listed in Table \ref{param}.

\begin{figure}
\centering
	\subfigure[Parameter estimation by using $v-i$ characteristics]{
		\includegraphics[scale=0.56]{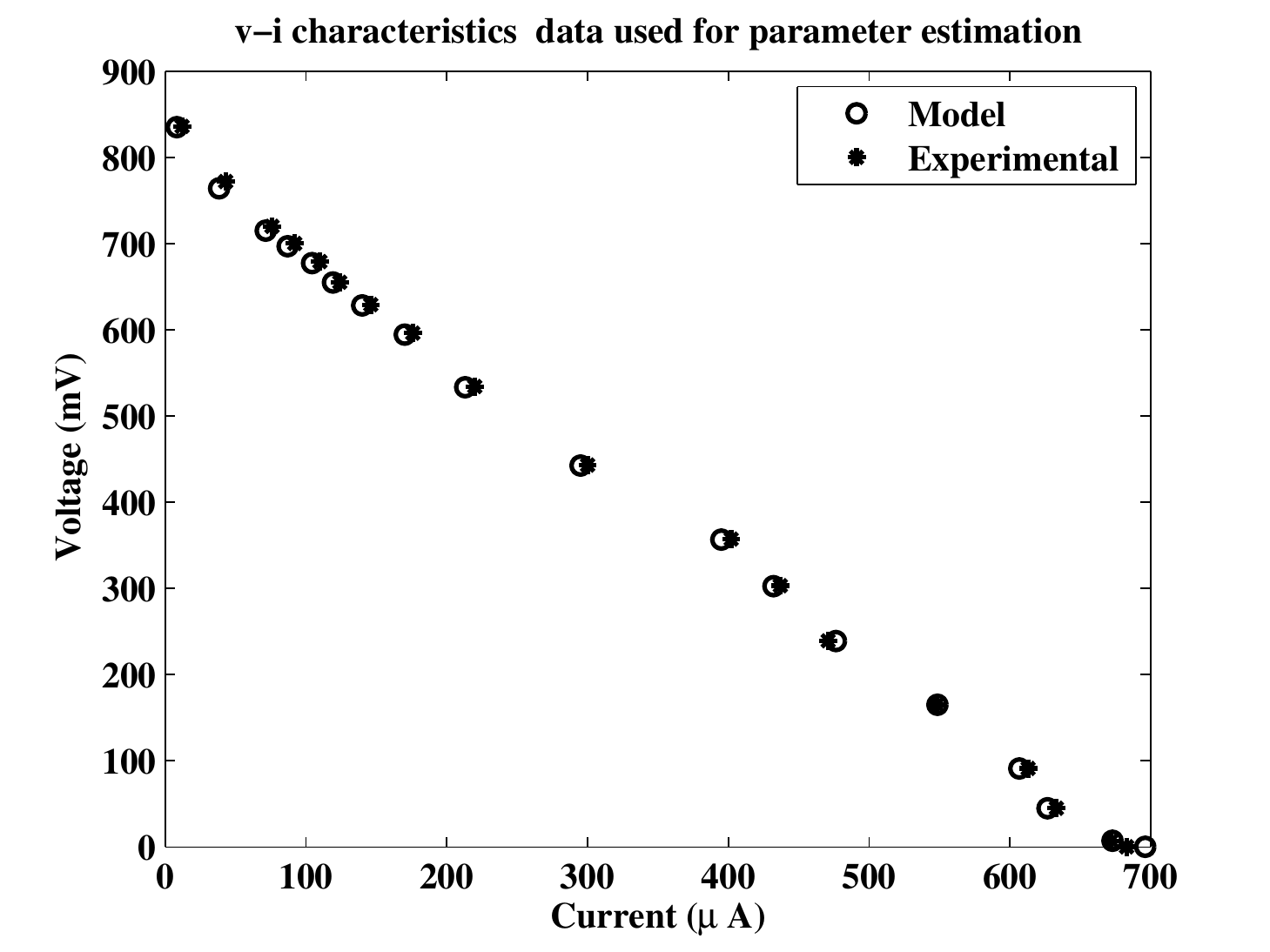}\label{vi_chara}}
	\subfigure[Power density curve for the points used for parameter estimation $\mu PSC$ ]{
		\includegraphics[scale=0.56]{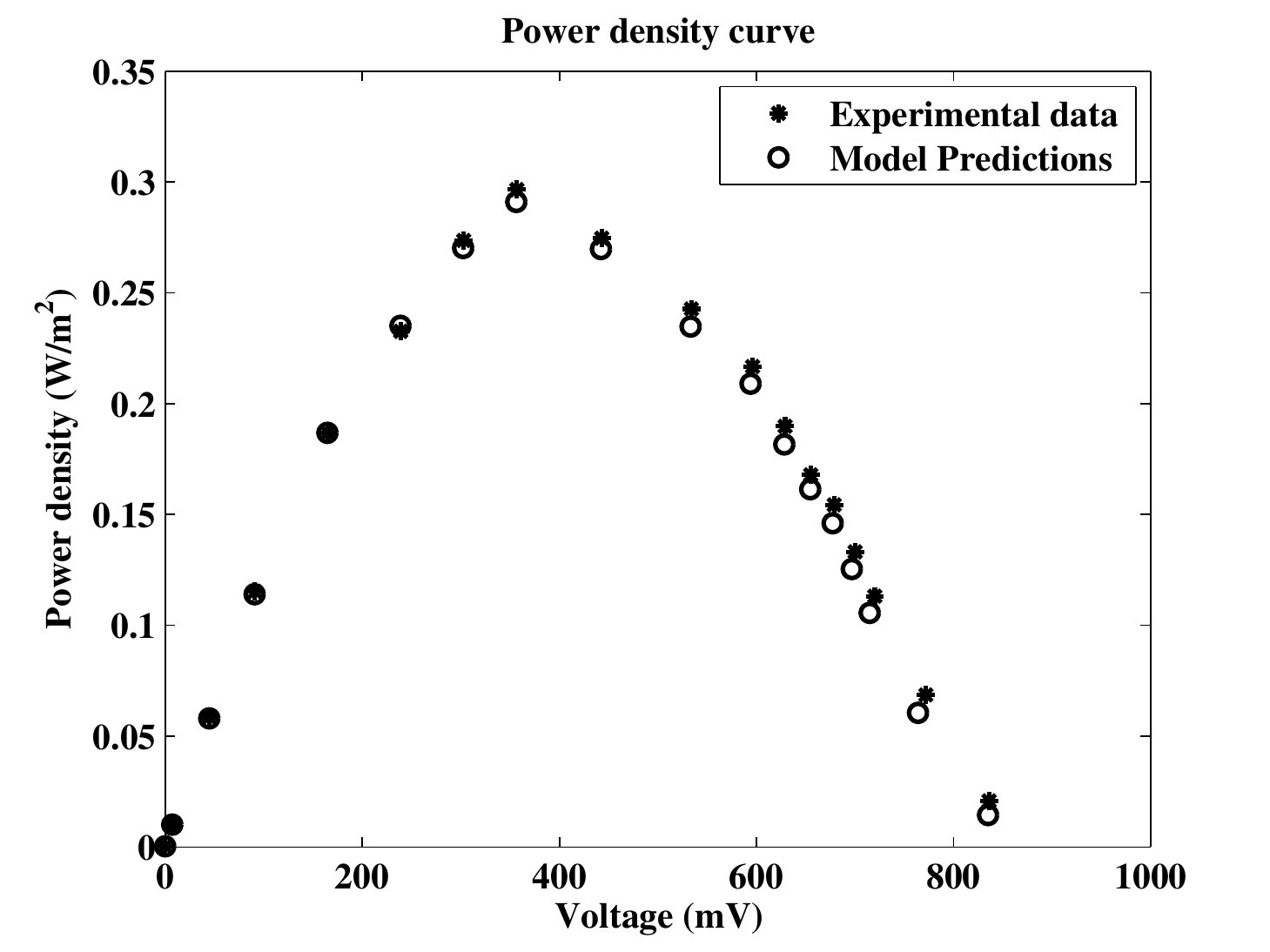}\label{powerdens}}
	\caption{Parameter estimation. \Cref{vi_chara} represents the predicted v-i characteristics by using the estimated parameters and \Cref{powerdens} represents corresponding power densities compared to experimental data }
\end{figure}

\begin{center}
\begin{table*}[!htbp]
\caption{Parameter values used for simulation of $\mu PSC$} \label{param}
\begin{tabular}{p{0.12\linewidth} p{0.35\linewidth}p{0.14\linewidth}p{0.14\linewidth}p{0.14\linewidth}}
\hline
Parameter &  Description  &  Value  &  Units  & Source  			\\
\hline
$x_{0}$	& Initial concentration of cells  &  12.2  &$\frac{g}{m^{3}} $& \cite{arvind2014advancedfabric}\\[3ex]
$k_{2}$	& Cell death rate constant  &  $5.32\times 10^{-6}$  &$\frac{1}{sec} $&  \cite{thornton2010modeling}\\[3ex]
$N_{0}$	&  Initial concentration of nutrients	&  2890  &$\frac{g}{m^{3}} $&\cite{arvind2014advancedfabric} \\[3ex]
${Y_{x,N}}$ &  Yield coefficient of cells w.r.t nutrients  & 10 & Dimensionless & \cite{thornton2010modeling}\\[3ex]
$\mu_{max}$ &  Maximum growth rate of C.reinhardtii  & $5 \times 10^{-5} $ & $sec^{-1}$ & \cite{vitova2011chlamydomonas}\\[3ex]
$K_{N}$ & Half saturation constant of cell growth w.r.t nutrients   & {4} &  $\frac{g-nutrient}{m^{3}}$   & \cite{thornton2010modeling}\\[4ex]
$R_{int} $ & Internal resistance & 599 & Ohm & Calculated  from \cite{arvind2014advancedfabric} \\ [3ex]
$R$ & Universal gas constant & 8.314 & $\frac{g}{mol-K}$ & \\ [3ex]
$T$ & Temperature & 298 & $K$ & \\[3ex]
$F$ & Faraday constant & 96486 &  $\frac{col}{mol}$ &  \\ [3ex]
$N_{Av}$ & Avogadro number & $6.023 \times 10^{23}$ &  $\frac{number of paricles}{mol of species}$ & \\ [3ex]
$L_0$ & Light intensity  & 625  & lux & \cite{arvind2014advancedfabric} \\[3ex] 
$C_f$ & Conversion factor & $1\times 10^{16}$ & $\frac{photons}{sec}$ &\cite{convfact} \\[3ex] 
$Q$ & Quantum yield & 0.742& $\frac{electrons}{photons}$ &\cite{arvind2014advancedfabric} \\[3ex] 
$A_s$ & Light irradiation surface area& $6.25\times 10^{-4}$& $m^2$ &\cite{arvind2014advancedfabric} \\[3ex] 
$A_E$ & Electrode surface area & $4.84\times 10^{-4}$& $m^2$ &\cite{arvind2014advancedfabric} \\[3ex] 
$\eta_{eff}$ & Efficiency of uptake of photons by cells & 0.5& Dimensionless & Assumed \\[3ex] 
$x_{max}$ & Maximum cell density & $1\times 10^{5}$& $\frac{g}{m^3}$ & \cite{thornton2010modeling} \\[3ex]
$\alpha$ & Parameter & 0.005& Dimensionless & Optimized \\[3ex]
$E_{cell}^0$ & Standard thermodynamic voltage & 1.241& $V$ & Calculated \\[3ex]
$n$ & Number of electrons transferred & 2& Dimensionless & Assumed \\[3ex]
\hline
\end{tabular}
\end{table*}
\end{center}
An optimization problem is solved to find the optimal values of $K$ for a chosen set of points in the $v-i$ data.
The available $v-i$ data ($32$ points) is divided in to two sets, first set($18$ points) is used for parameter estimation and the second set($14$ points), test data, is used to validate the model with the chosen parameters.  
The objective of the optimization problem is to reduce the sum square error (SSE) between the experimental and the predicted $v-i$ values from model. The optimization problem formulated is shown in \eqrefn{optiprob}.
\begin{multline}\label{optiprob}
\underset{K(i)}{\text{minimize}} 
\hspace{1cm}Obj = \sum_i(E_{exp}-E_{model}(K(i)))^2 +\\ \qquad \sum_i(I_{exp}- I_{model}(K(i)))^2\\
\text{Subject to}\, K(i) \geq 0,~i=1,2:2:30,31,32. 
\end{multline}
Optimum values of characteristic rate constant $K$  are estimated for the $18$ chosen points from the $v-i$ data. Since the model developed is a lumped parameter model, the effects of phenomena that are not modeled have been incorporated through the variation of the characteristic rate constant for each external load. Figure \ref{vi_chara} shows the comparison between the experimental data points of  $v-i$ data used for parameter estimation and the $v-i$ values obtained from the model with the estimated parameters. The RMSE of the fit is 0.0025. This indicates that the parameters estimated are consistent and accurate.

Power densities are calculated based on the $v-i$ data from the model and are plotted against the experimental power densities. Figure \ref{powerdens} compares power density calculated from model and experimental values. The fit emphasizes the capability of the model to produce consistent output with the trained data and the optimized parameters.   

\begin{figure}
\centering
\includegraphics[scale=0.56]{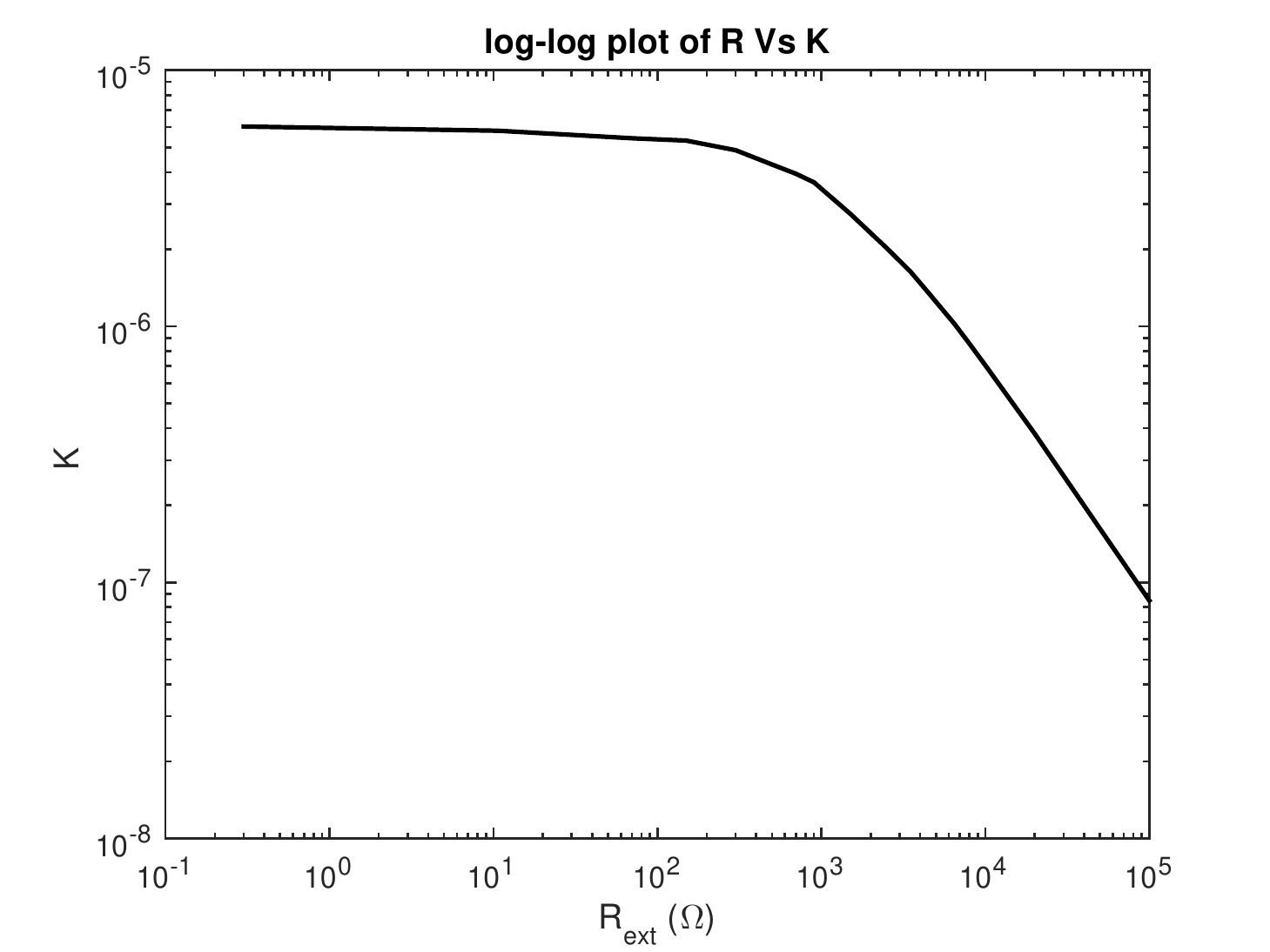}
\caption{Log-Log plot of $K$ vs $R_{ext}$}
\label{loglogka}
\end{figure}

Figure \ref{loglogka} shows the log-log plot of the estimated $K$ and $R_{ext}$. It is interesting to observe that there are two power law regions in the plot. The implications of this behavior of $K$ as a function of $R_{ext}$ provide some insights about the operating regimes of the cell. A discussion on these insights are presented later.
   
\begin{figure}
\centering
	\includegraphics[scale=0.56]{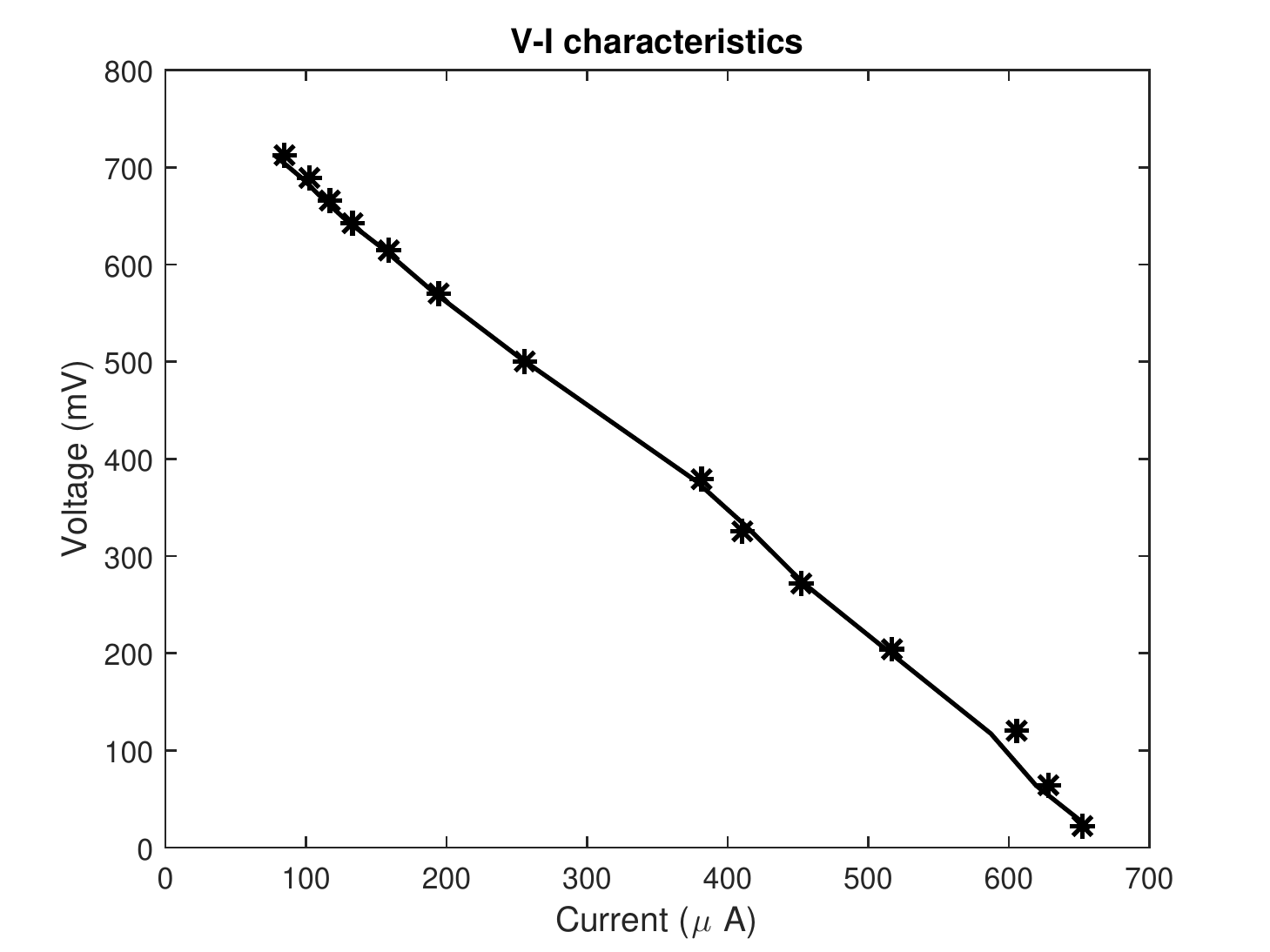}
	\caption{v-i characteristics prediction for test data}
	\label{vi_test}
\end{figure}
\section{Model validation}\label{modelvalid}

Model validation is a crucial step in any modeling exercise. In the present work the model is validated by using the steady-state experimental $v-i$ data, which is not used in parameter estimation. 

For $v-i$ data validation, the model responses are predicted for the $14$ test data points.  Figure \ref{vi_test} shows comparison between the experimental and the predicted $v-i$ characteristics. The RMSE is 0.0024 for the fit. The fit suggests that the model is able to quite accurately predict the voltage and current for the test data. The K values for the test data are estimated as a non-parametric interpolation of the K data identified from the training data. 

\section{Sensitivity analysis} \label{sensitivity}
\begin{figure}
\centering
\subfigure[$v-i$ characteristics for different electrode surface areas of  $\mu PSC$.]{\includegraphics[scale =0.51]{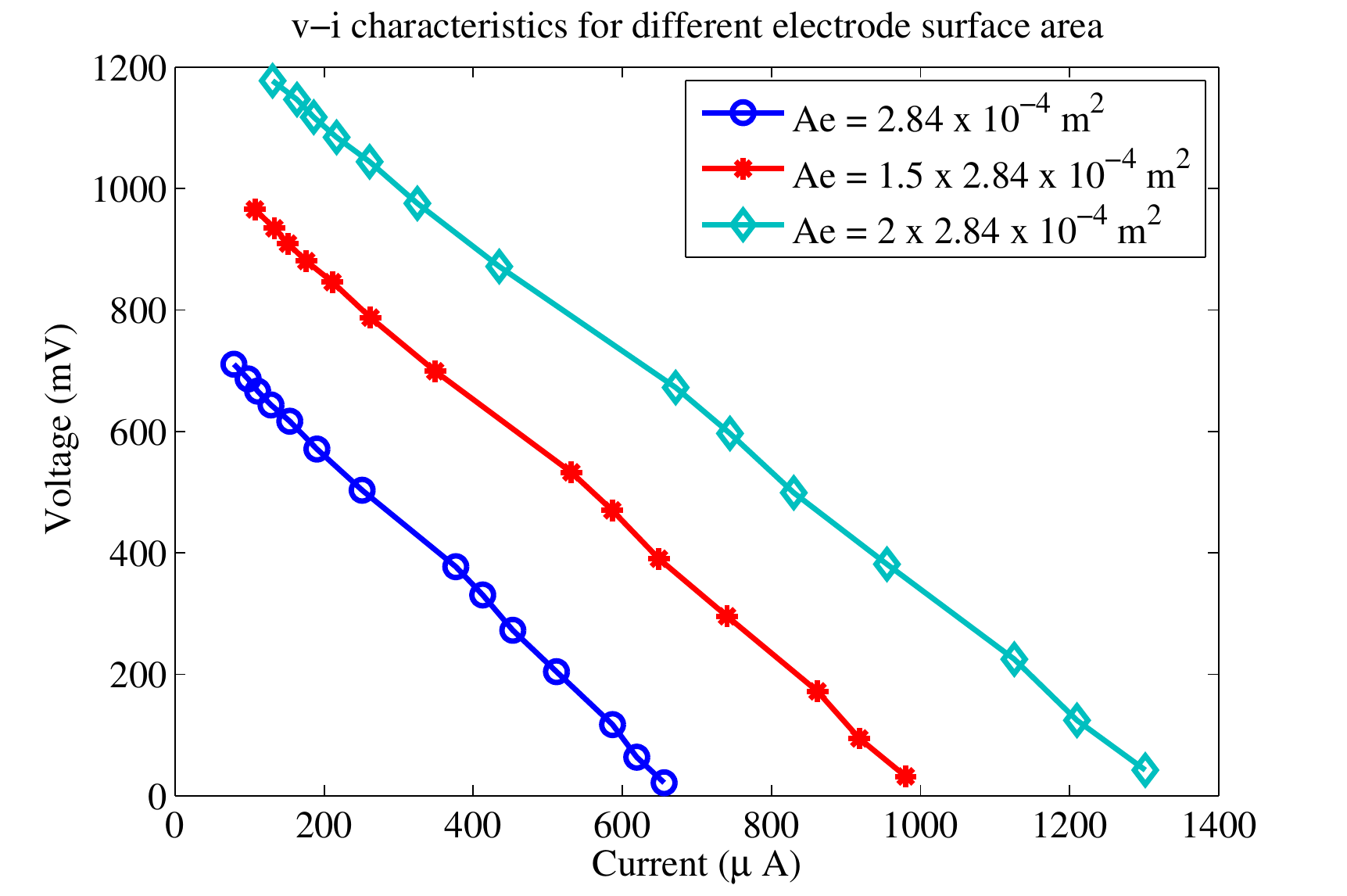}\label{Ae_VI}}
\subfigure[$v-i$ characteristics for different incident light intensities in $\mu PSC$.]{\includegraphics[scale =0.51]{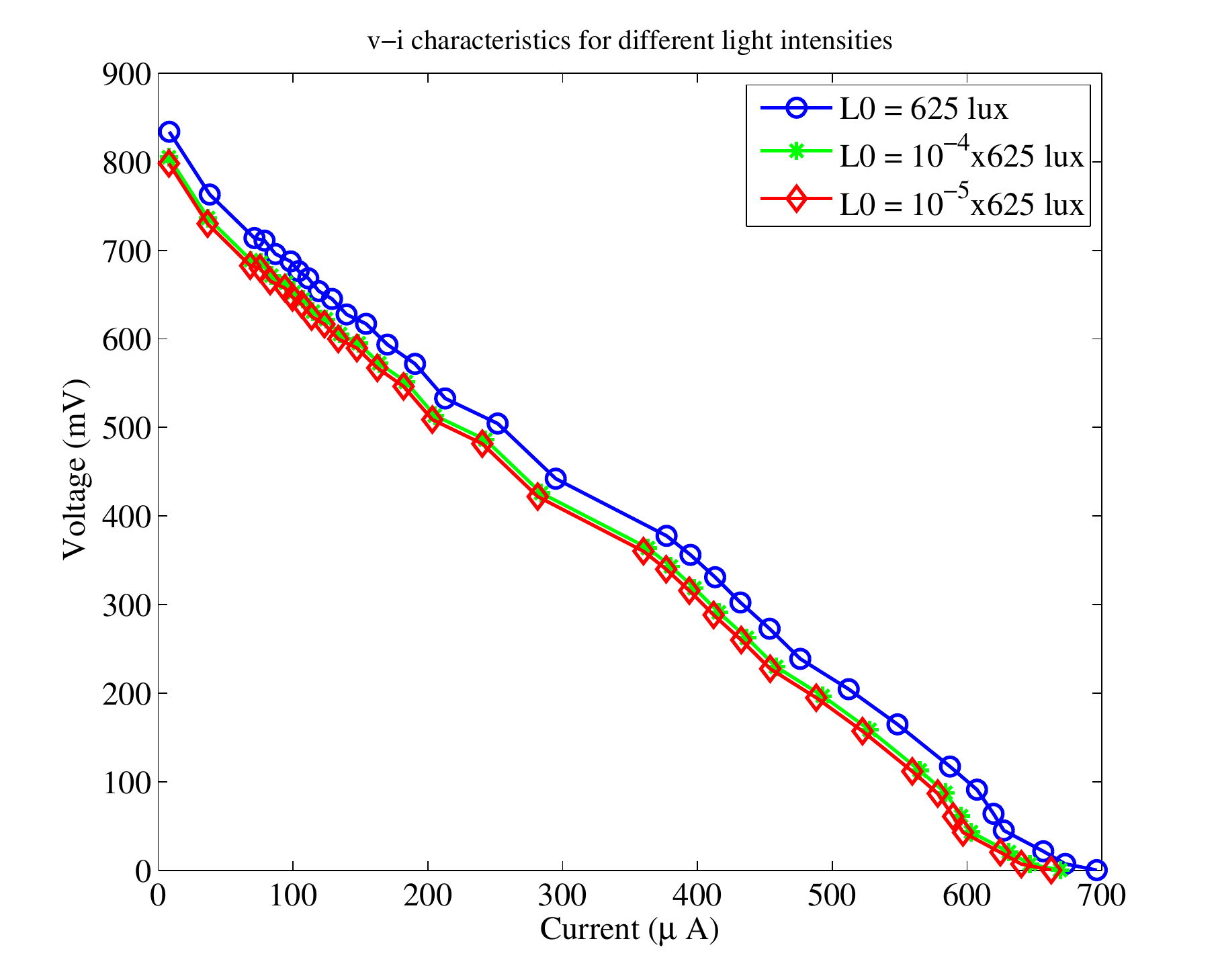}\label{L0_VI}}
\caption{Sensitivity analysis of $\mu PSC$ by varying electrode surface area ($A_e$) and light intensities ($L_0$)}
\end{figure}
Sensitivity analysis is performed to understand the opportunities for performance assessment. The proposed model is used for predicting the performance of the device for changes in:
\begin{inparaenum}[\itshape a\upshape)]
\item Electrode surface area $(A_e)$,
\item Incident light intensity$(L_0)$,
\item the concentration of cells in anode chamber$(x_0)$,
\item Light exposure surface area $(A_s)$.
\end{inparaenum} \label{applications}
The model predicted $v-i$ characteristics for  different electrode surface areas are shown in the Figure \ref{Ae_VI}. The model predicts increase in current with increase in electrode surface area which is consistent. A similar trend was also obtained for power density. Decrease in electrode surface decreases the area available for reactions to occur directly affecting the performance. This shows that even with the same cell concentrations much better performance might be possible by increasing the electrode surface area.

\begin{figure}
\centering
\subfigure[$v-i$ characteristics for different illumination surface areas of  $\mu PSC$]{\includegraphics[scale =0.51]{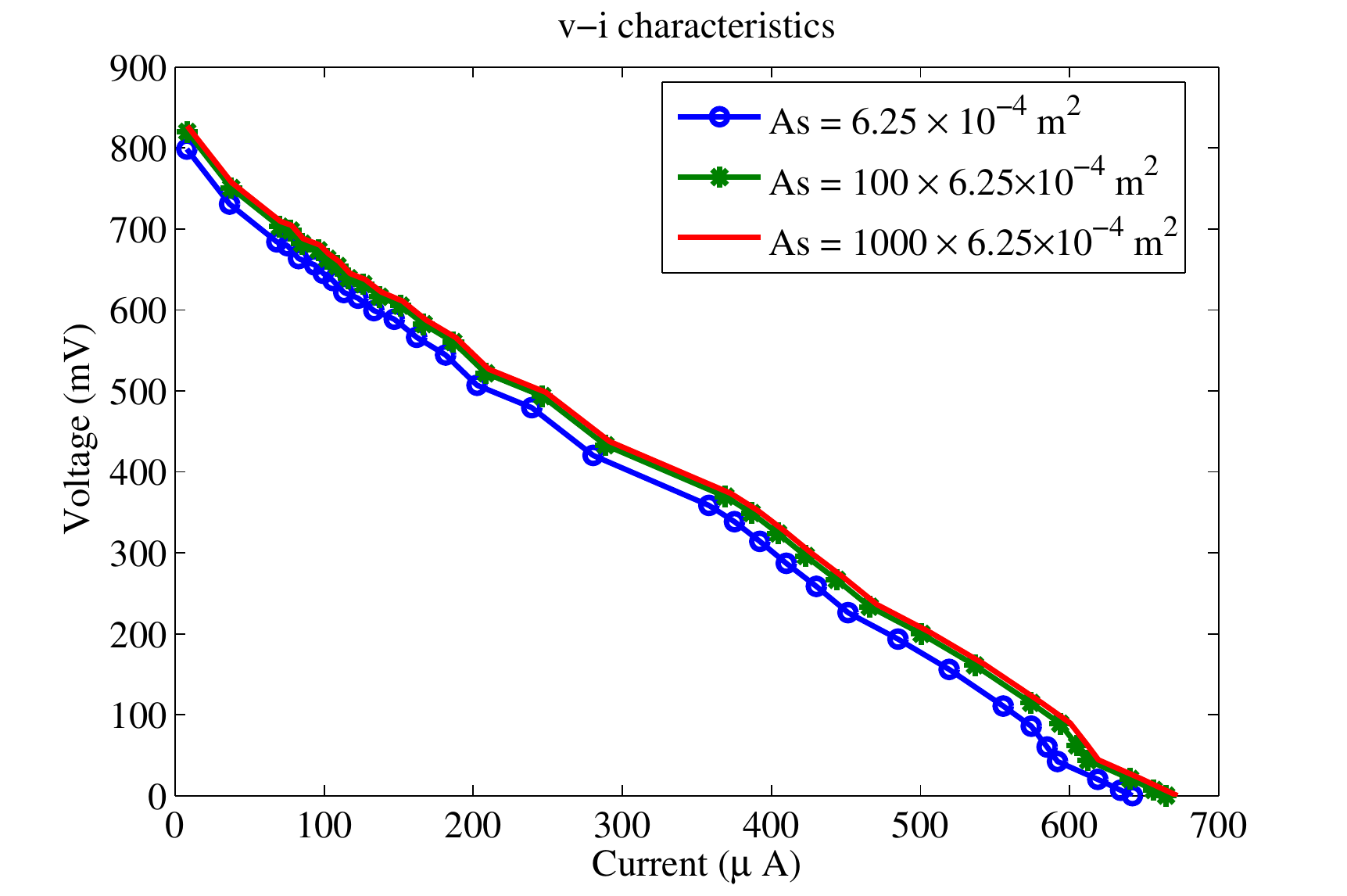}\label{As_VI}}
\subfigure[$v-i$ characteristics for different initial cell concentrations of $\mu PSC$]{\includegraphics[scale =0.5]{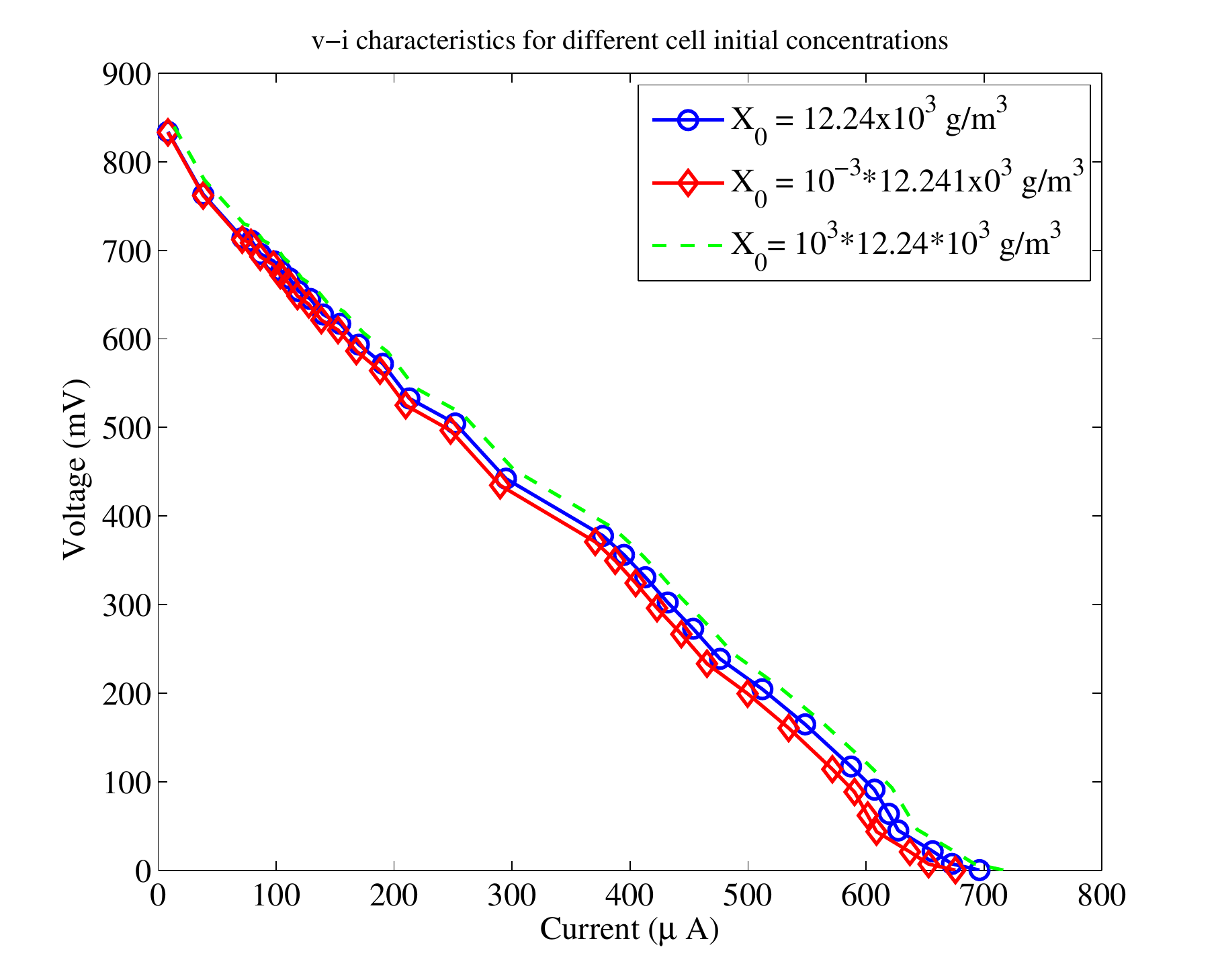}\label{x0_VI}}
\caption{Sensitivity analysis of $\mu PSC$ by varying initial cell concentrations ($X_0$) and illumination surface area ($A_s$)}
\end{figure}

The $v-i$ characteristics for the response of the system to different light intensities is shown in \Cref{L0_VI}. The model predicts decrease in current with decrease in light intensity and vice versa. One would expect that the current produced should be a very strong function of $L_0$, but the results show that the current is a weak function of the former. In other words, for this  $\mu PSC$, between increasing the electrode area and light intensity (or illumination surface area as shown in Figure \ref{As_VI}), the former is preferable. \Cref{x0_VI} shows the $v-i$ characteristics for various initial cell concentrations. The decrease in the  current and voltage is observed  with decreasing cell concentrations as shown in  \Cref{x0_VI}. It is observed that the model predicts the  current as a weak function of cell concentration. This is consistent with the other results as the strong dependence on the electrode area shows that the current cell population itself is underutilized and increasing the cell count is not likely to increase current because of the paucity of the active surface area for the electrochemical reactions. It is important to note that these conclusions are not easy to reach by looking at just the $v-i$ data without this modeling exercise. This emphasizes the power of this simple model in identifying key performance limiting factors.
\section{Discussions} \label{discuss}
\begin{figure}
\centering
	\subfigure[Log-Log plot of $K$ and $R_{ext}$ of $\mu PSC$ for $\alpha = 0.005$]{
		\includegraphics[scale=0.5]{loglog_vi}\label{goodka}}
	\subfigure[Log-Log plot of $K$ and $R_{ext}$ for $\alpha = 0.5$ showing discontinuities in parameter estimation]{
		\includegraphics[scale =0.5]{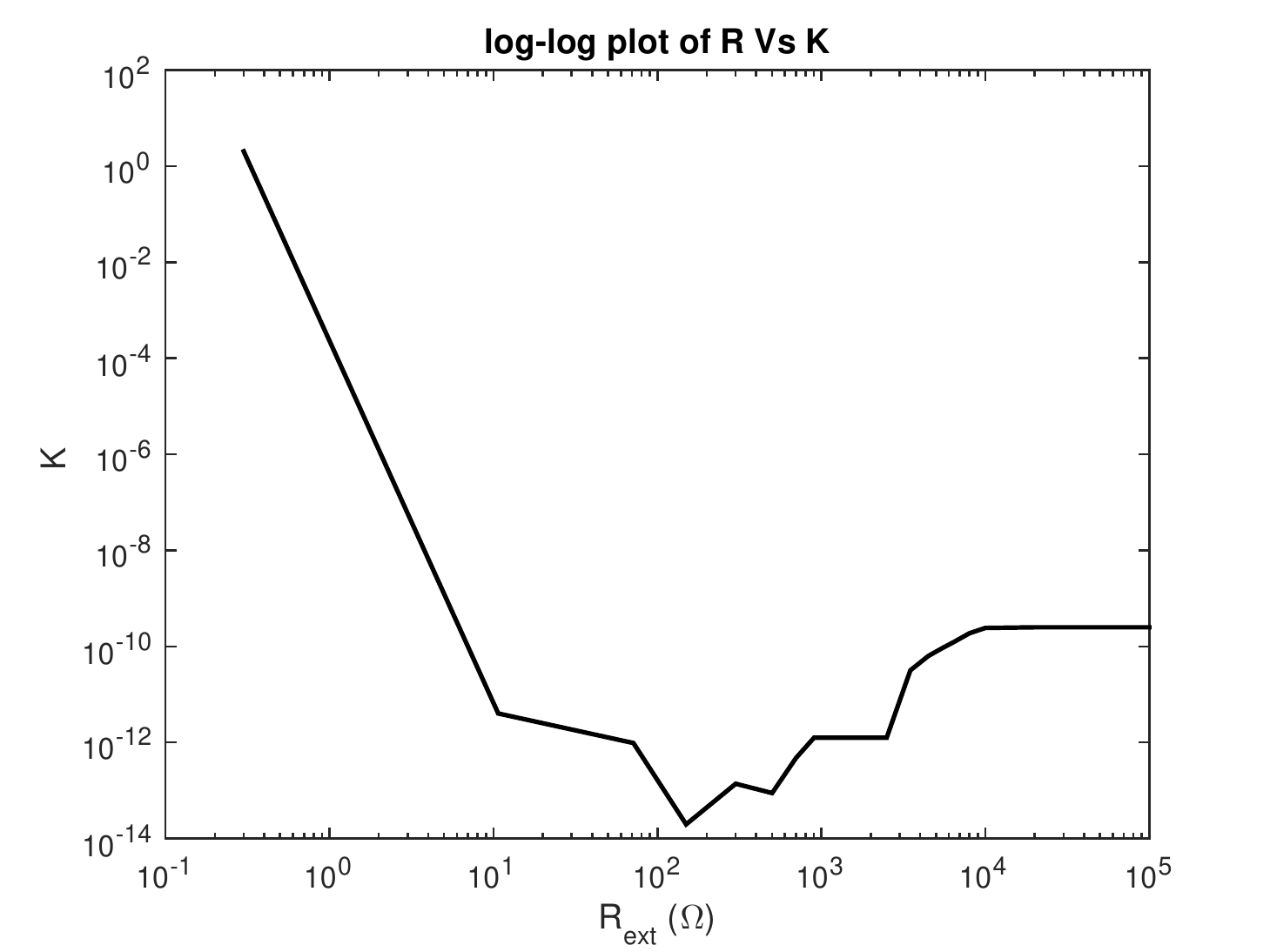}\label{worstk}}
	\caption{Comparison of the estimated parameter K at two different values of $\alpha$.}
\end{figure}

When $K$ is plotted vs $R_{ext}$ on a log-log graph (see Figure \ref{goodka}.) it can be observed that the experimental data collected for the cell can be divided into two regions. In the first region, the order of $K$ is almost constant ,i.e, the cell is being operated in an ohmic region, where the performance is not mass transfer limited. The reactants are supplied at electrodes at the same rate at which the reactants are used up in reaction. The second regime where  $K$ decreases as a power law, corresponds to the reaction rate limiting region (activation regime) of the operation of the cell. The sudden drop of voltage from $E_{0}$ in v-i profile, which corresponds to the activation loss dominant regime, is captured by large changes in parameter K. 

The key parameter in the $\mu PSC$ model is K, which represents the several transport phenomena and the rates of reactions. The K is related to $\alpha $ through $K=k_a^\alpha~k_c^{(1-\alpha)}$. $\alpha$ also has an interpretation of a charge transfer coefficient in Butler-Volmer equation, which is typically assumed to be $0.5$. However, in our optimization approach whenever the range of $\alpha$ was fixed as $O(1)$, there were discontinuities  in the estimated values of $K$ as shown in the log-log plot of $K$ vs $R_{ext}$ in Figure \ref{worstk}. However, for a small $\alpha$ value as used in our model, we could observe the natural two power law region curves as shown in Figure \ref{vi_test}. Hence, a value of $\alpha=0.005$  was chosen, which is still in the acceptable range of $0-1$. It is also worthwhile to keep in mind that the complicated multi-step photosynthesis reaction mechanism has been simplified into a single water splitting reaction step and hence the charge transfer coefficient can only be thought of as another parameter in the model.

It is also interesting to observe that when $\alpha = 0.5$ is used to fit the v-i data, K is almost of the same order of magnitude in the lower current regime and then varies significantly in the high current regime. That is, $\alpha = 0.5$ fits the data into the regimes where ohmic losses and mass transfer losses are dominant. However, for $\alpha = 0.005$, data was fit to regimes where ohmic losses and activation losses are dominant as mentioned before. The reason for a smooth fit for $\alpha = 0.005$ compared to $\alpha = 0.5$ could be because the experimental data that was used for parameter estimation lies in the region of dominant ohmic and activation losses.

\section{Conclusions} \label{conclude}
A mathematical model to predict the performance of a  $\mu PSC$ was developed. The model was thoroughly validated with steady-state experimental $v-i$ data. It was shown that this model could be used to predict the behavior of the $\mu$PSC that was considered. Several insights provided by the model regarding the operation of the  $\mu PSC$ were described. For this particular design, the model was able to unequivocally identify that the performance limitation as largely due to lack of enough active sites for reaction and not cell concentration or light intensity. As future work, if the model is extended to include the geometry of the electrode patterns, diffusion phenomena, and the multi-step reaction processes that occur, then it can be used to comprehensively optimize the various design and operational parameters of a $\mu$PSC.  

\section*{Acknowledgment}
The authors would like to thank the PEER group, Dr. Simona Badilescu and Dr.~Jayan Ozhi Kandthil from Concordia University, Montreal, Canada for their support. The authors also thank Mr.~Laya Das from IIT Gandhinagar, and Mr.~Srinivasan Raman and Dr.~Parham Mobed from Texas Tech university for their help in discussions on solutions to model equations.


\begin{thebibliography}{10}
	\providecommand{\url}[1]{#1}
	\csname url@samestyle\endcsname
	\providecommand{\newblock}{\relax}
	\providecommand{\bibinfo}[2]{#2}
	\providecommand{\BIBentrySTDinterwordspacing}{\spaceskip=0pt\relax}
	\providecommand{\BIBentryALTinterwordstretchfactor}{4}
	\providecommand{\BIBentryALTinterwordspacing}{\spaceskip=\fontdimen2\font plus
		\BIBentryALTinterwordstretchfactor\fontdimen3\font minus
		\fontdimen4\font\relax}
	\providecommand{\BIBforeignlanguage}[2]{{%
			\expandafter\ifx\csname l@#1\endcsname\relax
			\typeout{** WARNING: IEEEtran.bst: No hyphenation pattern has been}%
			\typeout{** loaded for the language `#1'. Using the pattern for}%
			\typeout{** the default language instead.}%
			\else
			\language=\csname l@#1\endcsname
			\fi
			#2}}
	\providecommand{\BIBdecl}{\relax}
	\BIBdecl
	
	\bibitem{siu2008microfabricated}
	C.-P. Siu and M.~Chiao, ``A microfabricated pdms microbial fuel cell,''
	\emph{Microelectromechanical Systems, Journal of}, vol.~17, no.~6, pp.
	1329--1341, 2008.
	
	\bibitem{cheng2011increasing}
	S.~Cheng and B.~E. Logan, ``Increasing power generation for scaling up
	single-chamber air cathode microbial fuel cells,'' \emph{Bioresource
		technology}, vol. 102, no.~6, pp. 4468--4473, 2011.
	
	\bibitem{zhang1995modelling}
	X.-C. Zhang and A.~Halme, ``Modelling of a microbial fuel cell process,''
	\emph{Biotechnology Letters}, vol.~17, no.~8, pp. 809--814, 1995.
	
	\bibitem{picioreanu2010modelling}
	C.~Picioreanu, K.~P. Katuri, M.~C. van Loosdrecht, I.~M. Head, and K.~Scott,
	``Modelling microbial fuel cells with suspended cells and added electron
	transfer mediator,'' \emph{Journal of applied electrochemistry}, vol.~40,
	no.~1, pp. 151--162, 2010.
	
	\bibitem{pinto2011unified}
	R.~P. Pinto, B.~Srinivasan, and B.~Tartakovsky, ``A unified model for
	electricity and hydrogen production in microbial electrochemical cells,'' in
	\emph{Preprints of the 18th Intenational Federation of Automatic Control
		(IFAC) world congress Milano (Italy) August}, 2011.
	
	\bibitem{yagishita1996photosynthetic}
	T.~Yagishita, S.~Sawayama, K.-i. Tsukahara, and T.~Ogi, ``Photosynthetic
	bio-fuel cells using cyanobacteria,'' \emph{Renewable energy}, vol.~9, no.~1,
	pp. 958--961, 1996.
	
	\bibitem{lam2004bio}
	K.~B. Lam, E.~Johnson, and L.~Lin, ``A bio-solar cell powered by sub-cellular
	plant photosystems,'' in \emph{Micro Electro Mechanical Systems, 2004. 17th
		IEEE International Conference on.(MEMS)}.\hskip 1em plus 0.5em minus
	0.4em\relax IEEE, 2004, pp. 220--223.
	
	\bibitem{rosenbaum2005utilizing}
	M.~Rosenbaum, U.~Schr{\"o}der, and F.~Scholz, ``Utilizing the green alga
	chlamydomonas reinhardtii for microbial electricity generation: a living
	solar cell,'' \emph{Applied microbiology and biotechnology}, vol.~68, no.~6,
	pp. 753--756, 2005.
	
	\bibitem{lam2006mems}
	K.~B. Lam, E.~A. Johnson, M.~Chiao, and L.~Lin, ``A mems photosynthetic
	electrochemical cell powered by subcellular plant photosystems,''
	\emph{Microelectromechanical Systems, Journal of}, vol.~15, no.~5, pp.
	1243--1250, 2006.
	
	\bibitem{chiao2006micromachined}
	M.~Chiao, K.~B. Lam, and L.~Lin, ``Micromachined microbial and photosynthetic
	fuel cells,'' \emph{Journal of Micromechanics and Microengineering}, vol.~16,
	no.~12, p. 2547, 2006.
	
	\bibitem{shahparnia2011polymer}
	M.~Shahparnia, ``Polymer micro photosynthetic power cell: Design, fabrication,
	parametric study and testing,'' Ph.D. dissertation, Concordia University,
	2011.
	
	\bibitem{arvind2014advancedfabric}
	A.~Ramanan, M.~Packirisamy, and S.~Williamson, ``Advanced fabrication,
	modeling, and testing of a micro-photosynthetic electrochemical cell for
	energy harvesting applications,'' \emph{Power Electronics, IEEE Transactions
		on}, vol.~PP, no.~99, pp. 1--1, 2014.
	
	\bibitem{ramanan2015advanced}
	A.~V. Ramanan, M.~Pakirisamy, and S.~S. Williamson, ``Advanced fabrication,
	modeling, and testing of a microphotosynthetic electrochemical cell for
	energy harvesting applications,'' \emph{Power Electronics, IEEE Transactions
		on}, vol.~30, no.~3, pp. 1275--1285, 2015.
	
	\bibitem{shahparnia2015micro}
	M.~Shahparnia, M.~Packirisamy, P.~Juneau, and V.~Zazubovich, ``Micro
	photosynthetic power cell for power generation from photosynthesis of
	algae,'' \emph{TECHNOLOGY}, vol.~3, no. 02n03, pp. 119--126, 2015.
	
	\bibitem{taiz2010plant}
	L.~Taiz and E.~Zeiger, \emph{Plant Physiology}.\hskip 1em plus 0.5em minus
	0.4em\relax Sinauer Associates, 2010.
	
	\bibitem{noren2005clarifying}
	D.~Noren and M.~Hoffman, ``Clarifying the butler--volmer equation and related
	approximations for calculating activation losses in solid oxide fuel cell
	models,'' \emph{Journal of Power Sources}, vol. 152, pp. 175--181, 2005.
	
	\bibitem{sunden2005transport}
	B.~Sund{\'e}n and M.~Faghri, \emph{Transport phenomena in fuel cells}.\hskip
	1em plus 0.5em minus 0.4em\relax WIT press, 2005, vol.~19.
	
	\bibitem{nelson2008lehninger}
	D.~L. Nelson, A.~L. Lehninger, and M.~M. Cox, \emph{Lehninger principles of
		biochemistry}.\hskip 1em plus 0.5em minus 0.4em\relax Macmillan, 2008.
	
	\bibitem{arvia2011electrochemical}
	A.~J. Arvia, A.~E. Bolzan, and M.~A. Pasquale, \emph{Electrochemical Catalysts:
		From Electrocatalysis to Bioelectrocatalysis}.\hskip 1em plus 0.5em minus
	0.4em\relax Wiley Online Library, 2011.
	
	\bibitem{gunawardena2008performance}
	A.~Gunawardena, S.~Fernando, and F.~To, ``Performance of a yeast-mediated
	biological fuel cell,'' \emph{International journal of molecular sciences},
	vol.~9, no.~10, pp. 1893--1907, 2008.
	
	\bibitem{thornton2010modeling}
	A.~Thornton, T.~Weinhart, O.~Bokhove, B.~Zhang, D.~M. Sar, K.~Kumar,
	M.~Pisarenco, M.~Rudnaya, V.~Savcenco, J.~Rademacher \emph{et~al.},
	``Modeling and optimization of algae growth,'' 2010.
	
	\bibitem{vitova2011chlamydomonas}
	M.~V{\'\i}tov{\'a}, K.~Bi{\v{s}}ov{\'a}, M.~Hlavov{\'a}, S.~Kawano,
	V.~Zachleder, and M.~{\v{C}}{\'\i}{\v{z}}kov{\'a}, ``Chlamydomonas
	reinhardtii: duration of its cell cycle and phases at growth rates affected
	by temperature,'' \emph{Planta}, vol. 234, no.~3, pp. 599--608, 2011.
	
	\bibitem{convfact}
	``Notes on practical photometry for image sensor and vision sensor
	developers,''
	http://www.ini.unizh.ch/~tobi/anaprose/recep/practicalPhotometry.pdf,
	accessed: 2010-11-13.
	
\end{thebibliography}
\end{document}